\documentclass[12pt,draft,leqno]{article}
\usepackage{amssymb, eucal, latexsym}

\newtheorem{theorem}{Theorem}[section]
\newtheorem{lm}[theorem]{Lemma}

\newtheorem{cor}[theorem]{Corollary}
\newtheorem{pro}[theorem]{Proposition}
\newtheorem{defi}[theorem]{Definition}

\newtheorem{nota}[theorem]{Notation}
\newtheorem{rem}[theorem]{Remark}

\newtheorem{fact}[theorem]{Fact}

\def\p{\varphi}
\def\a{\alpha}

\def\d{\delta}
\def\ep{\varepsilon}
\def\g{\gamma}
\def\GA{\Gamma}

\def\s{\sigma}
\def\SI{\Sigma}

\def\lra{\longrightarrow}

\def\sbe{\subseteq}

\def\stm{\setminus}
\def\ems{\emptyset}

\def\ex{\exists}
\def\fa{\forall}

\def\we{\wedge}

\def\ap{^{\,\prime}}
\def\inv{^{-1}}
\def\st{\ |\ }

\def\nin{\not\in}

\def\AA{{\cal A}}
\def\BB{{\cal B}}
\def\CC{{\cal C}}

\def\EE{{\cal E}}
\def\FF{{\cal F}}

\def\KK{{\cal K}}

\def\MM{{\cal M}}

\def\SS{{\cal S}}
\def\UU{{\cal U}}
\def\VV{{\cal V}}
\def\WW{{\cal W}}

\def\ZHC{{\bf Stone}}
\def\Bool{{\bf Bool}}

\def\2{\mbox{{\bf 2}}}
\def\3{\mbox{{\bf 3}}}

\def\int{\mbox{{\rm int}}}

\def\cl{\mbox{{\rm cl}}}

\def\doc{\hspace{-1cm}{\em Proof.}~~}
\def\sq{\hspace*{\fill} \hbox{\vrule\vbox{\hrule\phantom{o}\hrule}\vrule}}
\def\sqs{\sq \vspace{2mm}}

\def\NNNN{\mathbb{N}}
\def\DDDD{\mathbb{D}}
\def\QQQQ{\mathbb{Q}}

\def\IIII{\mathbb{I}}

\title{{\LARGE\bf
An internal topological characterization of the subspaces of Eberlein compacts and related compacts -- I}\\
\vspace{0.35cm}
{\large\bf Georgi D. Dimov}\thanks{This paper was supported by the
project no. 93/2011  of the
Sofia
University $``$St. Kl. Ohridski".
}\\
\vspace{0.25cm}
 {\footnotesize Dept. of Math. and
Informatics, Sofia University,  5 J. Bourchier Blvd., 1164 Sofia,
Bulgaria}
}

\author{}

\date{}

\begin{document}
\maketitle
\begin{abstract}
We obtain an internal topological characterization of the subspaces of Eberlein compacts (respectively, Corson compacts,  strong Eberlein compacts,
uniform Eberlein compacts, $n$-uniform Eberlein compacts).
\end{abstract}



\footnotetext[3]{{\footnotesize {\em E-mail address:}
gdimov@fmi.uni-sofia.bg}}


\baselineskip = \normalbaselineskip

\section{Introduction}


\bigskip

In 1977, A. V. Arhangel'ski\v{\i} proved that every metrizable space has a compactification which is an Eberlein compact
(see \cite[Theorem 14]{A}) (all necessary definitions are given in the next section).
 In 1982 (in a private communication), he posed the question:
$``$When does a space $X$ have a compactification $cX$ which is an Eberlein compact?". Since the closed subspaces of Eberlein compacts are Eberlein compacts,
 this question is equivalent to the following problem: find an internal characterization of  the subspaces of Eberlein compacts (note that H. P. Rosenthal \cite{Ro} gave an internal characterization of Eberlein compacts). In this paper we obtain such a characterization and, moreover, we characterize internally the subspaces of Corson compacts, of uniform Eberlein compacts, of $n$-uniform Eberlein compacts and of strong Eberlein compacts (see Theorems \ref{MainE} and \ref{MainSE}); all these results, with the exception of that about the subspaces of $n$-uniform Eberlein compacts, were announced (without any proofs) in my paper \cite{D1}. The cited above theorem of A. V. Arhangel'ski\v{\i}, as well as its recent generalizations obtained by T. Banakh and A. Leiderman \cite[Theorem 1]{BL} and B. A. Pasynkov \cite[Corollary 5]{P}, follow immediately from our results; the same is true for two assertions of J. Lindenstrauss \cite[Propositions 3.1 and 3.2]{L}. The paper contains also some results which were not announced in \cite{D1}; the main of them are: (1) the new characterizations of Eberlein compacts,  Corson compacts,  uniform Eberlein compacts,  $n$-uniform Eberlein compacts and strong Eberlein compacts (see Theorem \ref{unifthE} and Fact \ref{unifthSEC}), (2) the equalities $hl(X)=c(X)=hc(X)$ for subspaces $X$ of Eberlein compacts (see Corollary \ref{e6}), (3)  the characterization of the spaces which are co-absolute with (zero-dimensional) Eberlein compacts (see Theorem \ref{abseb}), (4) a new proof of the famous Ponomarev Theorem \cite{P1} giving a solution to the Birkhoff Problem 72 \cite{Bi}, and (5) a discussion of  the question whether each Tychonoff (or normal $T_1$) space with a uniform base (in the sense of P. S. Alexandrov \cite{Alex}) is homeomorphic to a subspace of some Eberlein compact (see Remark \ref{remunif}).
 The  notion of an {\em almost subbase}\/ (introduced in \cite{D1})   plays a central role in all these results and in the whole paper.

 Not all of the assertions announced in the paper \cite{D1} are proved here; the proofs of the remaining portion of them will be given in the second part of this paper.

We now fix the notation.

If $A$ is a set, we denote by $|A|$ its cardinality.
 If $(X,\tau)$ is a topological space and $M$ is a subset of $X$, we
denote by $\cl_{(X,\tau)}(M)$ (or simply by $\cl(M)$ or
$\cl_X(M)$) the closure of $M$ in $(X,\tau)$ and by
$\int_{(X,\tau)}(M)$ (or briefly by $\int(M)$ or $\int_X(M)$) the
interior of $M$ in $(X,\tau)$. The Alexandroff compactification of
a locally compact Hausdorff non-compact space $X$ is denoted
by $\a X$,
the  set of all positive natural numbers -- by  $\NNNN$,
  the real line (with its natural topology) -- by $\mathbb{R}$, and the subspaces $[0,1]$ and $\{0,1\}$ of $\mathbb{R}$ -- by $\IIII$ and $\DDDD$, respectively.
  As usual, $\omega=\NNNN\cup\{0\}$.

Let $X$ be a dense subspace of a space $Y$ and $U\sbe X$. The {\em extension of $U$ in $Y$}, denoted by $Ex_YU$, is the set $Y\stm\cl_Y(X\stm U)$; recall that if $U$ is open in $X$ then $Ex_YU$ is the greatest open subset of $Y$ whose trace on $X$ is $U$.

If  $X$ is a topological space and $f:X\lra \IIII$ is a continuous function, then we write, as usual, $Z(f)=f\inv(0)$ ({\em the zero-set of $f$}\/) and $coz(f)=X\stm Z(f)$ ({\em the cozero-set of $f$}\/). We set $Z(X)=\{Z(f)\st
f:X\lra \IIII$ is a continuous function$\}$ and $Coz(X)=\{coz(f)\st
f:X\lra \IIII$ is a continuous function$\}$.

We  denote by $\QQQQ_1$ the set of all positive rational numbers less than 1.

Recall that a subset $M$ of a topological space $X$ is
called {\em regular closed}\/ (respectively, {\em regular open}\/) if $M=\cl(\int (M))$ (respectively, $M=\int(\cl(M))$; we  denote by $RC(X)$ (respectively, $RO(X)$; $CO(X)$)
the collection  of all regular closed (respectively, regular open; clopen (= closed and open)) subsets of $X$.  Recall also that $RC(X)$ becomes a complete Boolean algebra
$(RC(X),0,1,\vee, \we, {}^*)$ under the following operations:
$ 1 = X$,  $0 = \emptyset$, $F^* = \cl(X\stm F)$, $F\vee G=F\cup G$,
$F\we G =\cl(\int(F\cap G)).$ Also, $RO(X)$ becomes a complete Boolean algebra
$(RO(X),0,1,\vee,\we,{}^*)$ under the following operations: $U\vee V=\int(\cl(U\cup V))$, $U\we V=U\cap V$,
$U^{*}=\int(X\setminus U)$, $0=\emptyset$, $1=X$. These two Boolean algebras are isomorphic.

The operation $``$complement" in Boolean algebras will be denoted
by $``$*".

We denote by
$S:\Bool\lra \ZHC$ the
Stone duality functor between the category $\Bool$ of Boolean algebras and Boolean
homomorphisms and the category $\ZHC$ of compact
zero-dimensional Hausdorff spaces  and
continuous maps (see, e.g. \cite{kop89}).

The class of all objects of a category $\CC$ is denoted by $|\CC|$.

All notions and notation which are not explained here can be found in \cite{E,Walker,kop89}.

{\bf By a $``$space" we  will mean a $``$topological $T_0$-space".}

\section{Preliminaries}
%

Let  $\GA$ be an index set and let $\IIII^\GA$ (respectively, $\DDDD^\GA$) be the Cartesian product of $|\GA|$
copies of $\IIII$ (respectively, $\DDDD$) with the Tychonoff topology. We set
$$\SI(\IIII,\GA) =\{ x\in \IIII^\GA\st |\{\g\in\GA\st x(\g) > 0\}|\le\aleph_0\},$$
$$\SI_*(\IIII,\GA)=\{ x\in \IIII^\GA\st (\fa \ep> 0) (|\{\g\in\GA\st x(\g)\ge\ep\}|<\aleph_0) \},$$
$$\s(\DDDD,\GA) =\{ x\in \DDDD^\GA\st |\{\g\in\GA\st x(\g)=1\}|<\aleph_0\},$$
and the topology on these subsets of $\IIII^\GA$ is the subspace topology. Obviously, $\s(\DDDD,\GA)\sbe\SI_*(\IIII,\GA)\sbe\SI(\IIII,\GA)$.

Let $\UU$ be a family of subsets of a set $X$. Recall that, for $n\in\omega$,  the {\em order}\/ of the family  $\UU$ is $\le n$ if any $n+2$ members of
the family $\UU$ have an empty intersection (i.e., each $x\in X$ is in at most $n+1$ members of $\UU$);
 the order of $\UU$ is {\em infinite}\/ if there is no $n\in\omega$ such that the order of $\UU$ is $\le n$.
 Let $n\in\NNNN$; the family $\UU$ is called {\em boundedly point-finite}\/ (respectively,  {\em $n$-boundedly point-finite}\/ (\cite{FMS})) if
 it has a finite order (respectively, is of order $\le n-1$). If $m\in\NNNN$ and the family $\UU$ is a union of countably many families $\UU_n$, where $n\in\NNNN$, we will say that $\UU$ is a {\em $\s$-point-finite}\/ (respectively, {\em $\s$-boundedly point-finite; $\s$-$m$-boundedly point-finite}\/) {\em family}\/  if all families $\UU_n$ are point-finite (respectively, boundedly point-finite; $m$-boundedly point-finite).
 The family $\UU$ is said to be {\em $T_0$-separating}\/ if, whenever $x\neq y$ are in $X$, then there exists $U\in\UU$ such that $|U\cap\{x,y\}|=1$; in this case we will also say  that the family $\UU$ {\em $T_0$-separates $X$}.
 Finally, when $X$ is a topological space,
the family $\UU$ is said to be {\em F-separating}
(see \cite{MR1}) if, whenever $x\neq y$ are in $X$, then there is some $U\in\UU$ such that $x\in U$ and $y\nin \cl_X(U)$, or vice versa.

 A compact Hausdorff space is called an {\em Eberlein
compact}\/ (briefly, EC), if it is homeomorphic to
a weakly compact  (i.e., compact in the weak topology) subset of a Banach space (\cite{L}).
D. Amir and J. Lindenstrauss proved that a compact
space is an Eberlein compact if and only if it can be embedded in
$\SI_*(\IIII,\GA)$ for some index set $\GA$  (see \cite[Theorem 1]{AL}). An internal characterization of Eberlein compacts was given by H. P. Rosenhthal \cite{Ro}.

\begin{theorem}\label{Rosenthal}{\rm (\cite{Ro})}
A compact Hausdorff space $X$ is an
Eberlein compact if and only if it has a $\s$-point-finite $T_0$-separating
collection of cozero-sets.
 \end{theorem}

A compact Hausdorff space that is homeomorphic to
a weakly compact subset of a Hilbert space is called a {\em uniform Eberlein compact}\/ (briefly, UEC); this class of spaces was introduced by Benyamini and Starbird in \cite{BS}. Let us recall the following result:

\begin{theorem}\label{argf}{\rm (\cite{BRW})}
Let $X$ be a compact Hausdorff space. Then  $X$ is a uniform Eberlein compact iff
$X$ has a $\s$-boundedly point-finite, $T_0$-separating family of cozero sets.
 \end{theorem}

Compact subspaces
of $\SI(\IIII,\GA)$ were called {\em Corson compacts}\/ (briefly, CC) by E.Michael and M.E.Rudin  \cite{MR1}. The internal characterization of Corson compacts is given below:

\begin{theorem}\label{mrcor}{\rm (\cite{MR1})}
Let $X$ be a compact Hausdorff space. Then  $X$ is a Corson compact iff
$X$ has a  point-countable, $T_0$-separating family of cozero sets.
 \end{theorem}

A compact space $X$ is called a {\em strong Eberlein compact}\/ (briefly, SEC) (\cite{Sim,BRW,Ro}) if it can be embedded in $\s(\DDDD,\GA)$ for some set $\GA$.
Equivalently, a compact space is a SEC iff it has a point-finite, $T_0$-separating family of clopen  sets.
K. Alster \cite{Al} proved that a space is a SEC iff it is a scattered EC.

We will need the following fundamental result:

\begin{theorem}\label{contim}{\rm (\cite{BRW,Gu,Sim,MR1})}
The classes of EC's, UEC's, SEC's and CC's are closed under continuous images.
 \end{theorem}

\begin{defi}\label{Urepr}{\rm (\cite{D1,D2})}
\rm
 Let $X$ be a space and $V$ a subset of it. If there
exists a collection $U (V) =\{U_n(V)\st n\in\NNNN\}$ such that $V =\bigcup\{U_n(V) \st
n\in\NNNN\}$ and $U_n(V)\sbe U_{n+1}(V)$, $U_{2n-1}(V)\in Z(X)$, $U_{2n}(V)\in Coz(X)$
for every $n\in\NNNN$, then we  say that the set $V$ is {\em U-representable}\/
and the collection $U(V)$ is a {\em U-representation of} $V$. If $\a$ is
a family of subsets of $X$ and, for every $V\in\a$, $U(V)$ is a
U-representation of $V$, then the family $U(\a)= \{U(V) \st V \in \a\}$
is called a {\em U-representation of} $\a$.
\end{defi}

The next lemma is standard:

\begin{lm}\label{Uf}{\rm (\cite{D2})}
 Let $X$ be a space and $U(V)=\{U_n(V)\st n\in\NNNN\}$ be a U-representa\-ti\-on of a subset $V$ of $X$. Then there exists a continuous function $f:X\lra\IIII$ such that
 $V=coz(f)$ and $f\inv([\frac{1}{2^{n-1}},1])=U_{2n-1}(V)$, for every $n\in\NNNN$.
\end{lm}

The definition given below was inspired by the results
of B. Efimov and G. \v{C}ertanov \cite{EC}:

\begin{defi}\label{alsubb}{\rm (\cite{D1})}
\rm
 A  family $\a$ of subsets of a space $X$ is said to be an {\em almost subbase
of} $X$ if there exists a U-representation $U(\a)$ of $\a$
 such that the family $\a\cup \{X\stm U_{2n-1}(V)\st V\in\a, n\in\NNNN\}$ is a subbase of $X$.
\end{defi}

\begin{rem}\label{rem0}
\rm
Every almost subbase $\a$ of a space $X$ consists of cozero-sets because a subset of $X$ is U-representable iff
it is a cozero-set. Obviously, a U-representable set $V$ can have many different
U-representations.

It is easy to see that a space has an almost subbase iff it is completely regular (see \cite{D1,D2}).

Let us note that in \cite{D2} the notion of almost subbase is introduced in a little bit different manner, namely, there is there an additional requirement that almost subbases have to be $T_0$-separating families. In \cite{D2} we work with arbitrary topological spaces. In this paper, as in \cite{D1}, we work with $T_0$-spaces only and this allows us to remove the extra requirement from \cite{D2}. Indeed, the condition that an almost subbase $\a$ of an arbitrary topological space $X$ is a $T_0$-separating family is used in \cite{D2} only  for showing that the family $\FF_\a=\{f_V\st V\in\a\}$, where $f_V$ is the function constructed in Lemma \ref{Uf} on the base of the given U-representation $U(V)$ of $V$,  separates points; if $X$ is a $T_0$-space, however, it is easy to see that the family $\FF_\a$ separates points even when $\a$ is not required to be a $T_0$-separating family (indeed, if $x\neq y$ are in $X$ and there is no $V\in\a$ such that $|V\cap\{x,y\}|=1$ then there is some $W\in\a$ and some $n\in\NNNN$ such that  $|(X\stm U_{2n-1}(W))\cap\{x,y\}|=1$; then $f_W(x)\neq f_W(y)$). Also, arguing in a similar way, we obtain that if $\a$ is an almost subbase of a $T_0$-space then the family $$\a\ap=\a\cup\{f_V\inv((r,1])\st V\in\a, r\in\QQQQ_1\}$$ (where the functions $f_V$ are constructed as above) is a $T_0$-separating family (it is even an F-separating family) and an almost subbase. Moreover, when $\a$ is a $\s$-point-finite (respectively, $\s$-locally finite; point-countable) family then $\a\ap$ has the same property.
\end{rem}

The next proposition (proved in \cite{D2}) generalizes the Nagata-Smirnov metrization theorem:

\begin{pro}\label{GNS}{\rm (\cite{D1,D2})}
 A space $X$ is metrizable iff it has a $\s$-locally finite almost subbase.
 \end{pro}

In \cite{D2}, using Proposition \ref{GNS},  the following theorem  was  proved as well:

\begin{theorem}\label{Baire}{\rm (\cite{D1,D2})}
Every Baire subspace of\/ $\SI_*(\IIII,\GA)$ contains a dense $G_\d$ metrizable subspace.
 \end{theorem}

\begin{cor}\label{gdelta}{\rm (\cite{BRW,Nam})}
 Every EC has a dense $G_\d$ metrizable subspace.
  \end{cor}



We will also need  the following simple lemma:


\begin{lm}\label{ecem}{\rm (\cite{EC})}
A continuous bijection $f : X \lra Y$ is a
homeomorphism if and only if there exists a subbase $\BB$ of $X$
such that $f(U)$ is open in $Y$ for every $U \in \BB$.
 \end{lm}

Let us recall the following results and definitions from \cite{D2}:

\begin{pro}\label{asbsp}{\rm (\cite{D2})}
Let $Y$ be a subspace of a space $X$ and let $\a$ be an almost subbase of $X$. Then $\a\cap Y=\{V\cap Y\st V\in\a\}$ is an almost subbase of $Y$.
 \end{pro}

 \begin{lm}\label{lem313}{\rm (\cite{D2})}
Let $X$ be a space and $\a$ be a $T_0$-separating family of
cozero-subsets of $X$. Let,  for
every $U\in\a$, a family   $\a_U=\{U_i\in Coz(X)\st i\in\NNNN\}$ is given,
such that  $U =\bigcup \{U_i \st i\in\NNNN\}$ and $U_i\sbe \cl(U_i)\sbe U_{i+1}$   for every $i\in\NNNN$. Then $\a^*=\bigcup\{\a_U\st U\in\a\}$ is an F-separating family of
cozero-subsets of $X$.
 \end{lm}

\begin{defi}\label{unifasb}{\rm (\cite{D2,D3})}
\rm
 A  family $\a$ of cozero-subsets of a space $X$ is said to be a {\em uniform almost subbase
of}\/ $X$ if, for every U-representation $U(\a)$ of $\a$,
  the family $\a\cup \{X\stm U_{2n-1}(V)\st V\in\a, n\in\NNNN\}$ is a subbase of $X$.
\end{defi}

\begin{lm}\label{lem315}{\rm (\cite{D2})}
If $X$ is compact then every F-separating family $\a$ of
cozero-subsets of $X$ is a uniform almost subbase of $X$.
 \end{lm}

\begin{theorem}\label{unifth}{\rm (\cite{D2,D3})}
A compact Hausdorff space $X$ is an Eberlein compact iff $X$ has a $\s$-point-finite (uniform) almost subbase.
 \end{theorem}

The next notion was introduced recently by B. A. Pasynkov \cite{P}:

\begin{defi}\label{nunif}{\rm (\cite{P})}
\rm
Let $n\in\omega$. A compact Hausdorff space $X$ is called {\em $n$-uniform Eberlein}\/ (briefly, {\em $n$-UEC}\/) if $X$ has  a $T_0$-separating family $\g=\bigcup_{i\in\NNNN}\g_i$ of cozero sets  such that the order of all families $\g_i$, for $i\in\NNNN$, is $\le n$ (i.e., in our terms, if $X$ has a $T_0$-separating, $\s$-$(n+1)$-boundedly point-finite family of cozero-sets).
\end{defi}

B. A. Pasynkov announced in \cite[Final remark]{P} that for any $n\in\NNNN$, there exists an $n$-UEC that is not an $(n-1)$-UEC.

\section{On Eberlein spaces and related spaces}

\begin{defi}\label{strpf}{\rm (\cite{Pe1,EP,D1})}
\rm
A family $\a$ of subsets of a set $X$ is called {\em strongly point-finite}\/ (respectively, {\em strongly point-countable}\/) if every subfamily $\mu$ of $\a$ with $|\mu|=\aleph_0$ (respectively, $|\mu|=\aleph_1$) contains
a finite subfamily $\mu\ap$ with empty intersection.

A family of subsets of a set $X$ is called  {\em $\s$-strongly point-finite}\/ if it is a union of countably many  strongly point-finite families.
\end{defi}

\begin{rem}\label{rem1}
\rm
 Obviously, a family $\a$ of subsets of a set $X$ is strongly point-finite iff it contains no infinite  subfamily with finite intersection property. In this form, the notion of  $``$strongly point-finite family" is introduced
  independently in \cite{Pe1} and \cite{EP} (in \cite{Pe1} there is no name for such families and in \cite{Pe2} they are called  $``$b-families");
  moreover,
  in \cite{Pe1}, this notion is
  attributed to A. V. Arhangel'ski\v{\i}.
  Of course, I was unaware of these facts when I introduced this notion in \cite{D1}.

   Clearly, a family $\a$ of subsets of a set $X$ is strongly point-countable iff it contains no uncountable  subfamily with finite intersection property. In this form, the notion of  $``$strongly point-countable family" is introduced in \cite{Pe1} (as $``$condition B1") but it is not studied there.

   Evidently, every strongly point-finite (respectively,  strongly point-countable) family is a point-finite (respectively,
 point-countable) family; the converse is not true.
  \end{rem}

\begin{lm}\label{Ex}
 Let $X$ be a dense subspace of a space $Y$, $n\in\omega$ and $\g$ be a strongly point-finite (respectively, boundedly point-finite, $n$-boundedly point-finite, strongly point-countable) family of open subsets of $X$. Then the family $$Ex_Y\g=\{Ex_YU\st U\in\g\}$$ is a strongly point-finite (respectively, boundedly point-finite, $n$-boundedly point-finite, strongly point-countable) family of open subsets of $Y$.
\end{lm}

\doc Let $\g$ be a strongly point-finite family of open subsets of $X$, $\mu\sbe Ex_Y\g$ and $|\mu|=\aleph_0$. Let $\mu_X=\{V\cap X\st V\in\mu\}$. Then $|\mu_X|=\aleph_0$ and $\mu_X\sbe \g$. Hence, there exists a finite subfamily $\mu_X\ap$ of $\mu_X$ with empty intersection. Since the operator $Ex_Y$ preserves  finite intersections, we get that $\mu\ap=Ex_Y\mu_X\ap$ is a finite subfamily of $\mu$ with empty intersection. Therefore, $Ex_Y\g$ is a strongly point-finite family of open subsets of $Y$.

The proof of the other three cases is analogous.
\sqs

In the proof of the case of Eberlein compacts of the next theorem, the constructions of the  family  $\a$  is taken from the Michael-Rudin proof of Rosenthal's theorem (see \cite[Theorem 1.4]{MR1}).

\begin{theorem}\label{unifthE}
Let $n\in\omega$ and $X$ be a compact Hausdorff space $X$.  Then $X$ is an EC (respectively, CC; UEC; $n$-UEC) iff
 $X$ has a $\s$-strongly point-finite (uniform) almost subbase  (respectively, strongly point-countable  (uniform)  almost subbase; $\s$-boundedly point-finite  (uniform)  almost subbase; $\s$-($n+1$)-boundedly point-finite  (uniform)  almost subbase).
 \end{theorem}

 \doc We  start with the proof for Eberlein compacts.

 \noindent($\Rightarrow$)  Let $X$ be an EC. By the Amir and Lindenstrauss theorem,
we can regard $X$ as a subspace of $\SI_*(\IIII,\GA)$ for some index set $\GA$.
Now,  for every $\g\in\GA$ and every $r\in\QQQQ_1$,
we set
$$V_{r,\g}=\{x\in X\st x(\g)> r\},\ \
 \a_r=\{V_{r,\g}\st \g\in\GA\},\ \ \a=\bigcup\{\a_r\st r\in\QQQQ_1\}.$$

 We will  prove that, for every $r\in\QQQQ_1$, $\a_r$ is a strongly point-finite family.

 Let $r\in\QQQQ_1$ and $\mu=\{V_{r,\g_n}\st n\in\NNNN\}$, where $\g_n\neq\g_m$ for $n,m\in\NNNN$ with $n\neq m$,  be an infinite countable subfamily of $\a_r$ consisting of non-empty sets.
 Set $\overline{\mu}=\{\cl_X(V_{r,\g_n})\st n\in\NNNN\}$.
 Then the family $\overline{\mu}$ has empty intersection. Indeed,
 since $\cl_X(V_{r,\g_n})\sbe\{y\in X\st y(\g_n)\ge r\}$ for every $n\in\NNNN$,
if some point $x=(x(\g))_{\g\in\GA}$ of $X$ belongs to every element of $\overline{\mu}$, then $x(\g_n)\ge r$ for every $n\in\NNNN$ and this is a contradiction because $r>0$
 and $X\sbe\SI_*(\IIII,\GA)$. So, $\bigcap\overline{\mu}=\ems$. Since $X$ is compact, there exists a finite subfamily $\overline{\mu_0}=\{\cl_X (V_{r,\g_{n_i}})\st i=1,\ldots,k\}$ (where $k$ is from $\NNNN$) of $\overline{\mu}$
 with empty intersection.
 Then the family $\mu_0=\{V_{r,\g_{n_i}}\st i=1,\ldots,k\}$ is a finite subfamily of $\mu$ with empty intersection.
 Therefore, $\a_r$ is a strongly point-finite family.

 It is easy to see that $\a$ is  an F-separating family. Then, by Lemma \ref{lem315}, $\a$ is a uniform almost subbase and, thus, it is an almost subbase.

\smallskip

\noindent($\Leftarrow$) Let $\a$ be a $\s$-strongly point-finite almost subbase of $X$. Then the family $\a\ap$ constructed in Remark \ref{rem0} is a $T_0$-separating, $\s$-strongly point-finite almost subbase of $X$. Since every strongly point-finite family is point-finite, our assertion follows from Theorem \ref{Rosenthal}.

\smallskip

We proceed with the proof for Corson compacts:

\smallskip

 \noindent($\Rightarrow$)
 Let $X$ be a CC. Then
we can regard $X$ as a subspace of $\SI(\IIII,\GA)$ for some index set $\GA$.
For every $\g\in\GA$ and every $r\in\QQQQ_1$, we define the set $V_{r,\g}$ as above, and
we set
$\a=\{V_{r,\g}\st \g\in\GA,\ r\in\QQQQ_1\}.$
We will  prove that $\a$ is a strongly point-countable family.
Let  $\mu=\{V_{r_\xi,\g_\xi}\st \xi<\omega_1\}$, where $(r_\xi,\g_\xi)\neq(r_\eta,\g_\eta)$ for $\xi,\eta<\omega_1$ with $\xi\neq \eta$,  be a subfamily of $\a$ consisting of non-empty sets. Set $\overline{\mu}=\{\cl_X(V_{r_\xi,\g_\xi})\st \xi<\omega_1\}$. Then the family $\overline{\mu}$ has empty intersection. Indeed, suppose that some point $x=(x(\g))_{\g\in\GA}$ of $X$ belongs to every element of $\overline{\mu}$.
 Since $\cl_X(V_{r_\xi,\g_\xi})\sbe\{y\in X\st y(\g_\xi)\ge r_\xi\}$ for every $\xi<\omega_1$,  we get that $x(\g_\xi)\ge r_\xi$ for every $\xi<\omega_1$. From the countability of $\QQQQ_1$, we obtain that there exists $\xi_0<\omega_1$ such that $|\{\xi<\omega_1\st r_\xi=r_{\xi_0}\}|=\aleph_1$. Hence $x(\g_\xi)\ge r_{\xi_0}>0$ for every $\xi<\omega_1$ such that $r_\xi=r_{\xi_0}$. This implies that $|\{\g\st x(\g)\neq 0\}|\ge\aleph_1>\aleph_0$. Since $x\in\SI(\IIII,\GA)$, we get a contradiction.
  Therefore, $\bigcap\overline{\mu}=\ems$. Now, the compactness of $X$ implies that there exists a finite subfamily
  $\overline{\mu}_0=\{\cl_X(V_{r_{\xi_i},\g_{\xi_i}})\st i=1,\ldots,k\}$ (where $k\in\NNNN$)  of $\overline{\mu}$ with empty intersection.
 Then the family $\mu_0=\{V_{r_{\xi_i},\g_{\xi_i}}\st i=1,\ldots,k\}$ is a finite subfamily of $\mu$ with empty intersection. Hence, $\a$ is a strongly point-countable family.

It is easy to see that $\a$ is  an F-separating family. Now, we finish the proof as in ECs' case.

\smallskip

\noindent($\Leftarrow$) It is completely analogous to the corresponding part of the proof for Eberlein compacts, however, Theorem \ref{mrcor} has to be used instead of
Theorem \ref{Rosenthal}.

\smallskip

The proof for uniform Eberlein  compacts is the following:

\smallskip

 \noindent($\Rightarrow$)
 Let $X$ be a UEC. Then $X$ has a  $\s$-boundedly point-finite, $T_0$-separating family $\a$ of cozero-subsets (see Theorem \ref{argf}). Let
 $\a^*$ be the family described in Lemma \ref{lem313}. Then $\a^*$ is an F-separating, $\s$-boundedly point-finite family. Thus, by Lemma \ref{lem315}, $\a^*$ is a $\s$-boundedly point-finite uniform almost subbase.

 \smallskip

\noindent($\Leftarrow$) It is completely analogous to the corresponding part of the proof for Eberlein compacts, however, Theorem \ref{argf} has to be used instead of
Theorem \ref{Rosenthal}.

\smallskip

The proof for $n$-uniform Eberlein  compacts is completely analogous to that for uniform Eberlein compacts, however, in it the definition of $n$-uniform Eberlein  compacts has to be used instead of Theorem \ref{argf}.
\sqs

\begin{defi}\label{espa}{\rm (\cite{D1})}
\rm
 Let $n\in\omega$. A space $X$ is called an {\em Eberlein space}\/ (brieily, {\em E-space}\/) (respectively, {\em C-space},
{\em UE-space}, {\em $n$-UE-space}\/) if $X$ has a $\s$-strongly point-finite (respectively, strongly point-countable, $\s$-boundedly
point-finite, $\s$-($n+1$)-boundedly point-finite) almost subbase. A space $X$ is said to be a {\em SE-space}\/ if it has a strongly
point-finite family $\a$ of clopen subsets such that $\a\cup \{X\stm V\st  V\in\a\}$ is a subbase of $X$.
\end{defi}

The assertions which are contained in the next theorem were announced in \cite{D1} with the exception of that about $n$-uniform Eberlein compactifications.

\begin{theorem}\label{MainE}{\rm (\cite{D1})}
 Let $n\in\omega$. A space X has a compactification $cX$ which is an Eberlein compact (respectively, Corson compact; uniform Eberlein compact; $n$-uniform Eberlein compact) iff it is an E-space
 (respectively, C-space; UE-space; $n$-UE-space).
\end{theorem}

\doc We  start with the proof for Eberlein compactifications.

\noindent($\Rightarrow$)  Let $X$ has a compactification $cX$ which is an EC. We can think that $X$ is a subspace of $cX$. By Theorem \ref{unifthE}, $cX$ has a $\s$-strongly point-finite almost subbase $\a\ap$. Then Proposition \ref{asbsp} implies that $\a=\a\ap\cap X$ is an almost subbase of $X$. Obviously, $\a$ is a $\s$-strongly point-finite family. Hence, $X$ is an E-space.

 \smallskip

 \noindent($\Leftarrow$) Let $X$ be an E-space. We will construct a compactification $cX$ of $X$ which is an Eberlein compact.

 Let $\a$ be a $\s$-strongly point-finite  almost subbase of $X$. For every $V\in\a$, we set $\IIII_V=\IIII$. Let $Z=\prod_{V\in\a}\IIII_V$. Let, for every  $V\in\a$, $f_V:X\lra I_V$ be the continuous function  corresponding to the given U-representation $U(V)$ of $V$ (see Lemma \ref{Uf}). Let $f=\Delta_{V\in\a}f_V:X\lra Z$ be the diagonal of the mappings $(f_V)_{V\in\a}$. By Remark \ref{rem0}, $f$ is a continuous injection. Let $cX=\cl_Z(f(X))$ and let, for any $V\in\a$, $\pi_V:Z\lra I_V$ be the projection.  For every $V\in\a$, define a function $f_V\ap:cX\lra I_V$  by the formula $f_V\ap=\pi_V|cX$.
 Then, for every $V\in\a$ and every $n\in\NNNN$, $f_V\ap(f(x))=f_V(x)$ for every $x\in X$, $f(V)=f(X)\cap coz(f_V\ap)$ and $f(X\stm U_{2n-1}(V))=f(X)\cap (f_V\ap)\inv([0,\frac{1}{2^{n-1}}))$;    hence, $f(V)$
  and $f(X\stm U_{2n-1}(V))$ are open subsets of $f(X)$.  Thus,
using Lemma \ref{ecem}, we get that
  $f$ is an embedding of $X$ into $cX$.  Therefore, $cX$ is a compactification of $X$. We will show that $cX$ is an Eberlein compact.

 Since, for any $V\in\a$, we have that $f(X)\cap coz(f_V\ap)=f(V)$,  we obtain that $coz(f_V\ap)\sbe Ex_{cX}(f(V))$. The family $\{f(V)\st V\in\a\}$ is $\s$-strongly point-finite. Thus, by Lemma \ref{Ex},
 the family $\{Ex_{cX}(f(V))\st V\in\a\}$ is $\s$-strongly point-finite. Therefore, the family $\a\ap=\{coz(f_V\ap)\st V\in\a\}$ is $\s$-point-finite. Adding to $\a\ap$ the family $\{(f_V\ap)\inv((r,1])\st V\in\a,\ r\in\QQQQ_1\}$, we get a $T_0$-separating, $\s$-point-finite family of cozero-subsets of $cX$. Hence, using Rosenthal's Theorem \ref{Rosenthal}, we obtain that $cX$ is an Eberlein compact.

 \medskip

 The proof for Corson (respectively, uniform Eberlein; $n$-uniform Eberlein) compactifications is completely analogous to that for Eberlein compactifications. However, in it we have to use Theorem \ref{mrcor} (respectively, Theorem \ref{argf}; the definition of $n$-uniform Eberlein compacts) instead of Theorem \ref{Rosenthal}.
 \sqs

\begin{cor}\label{subsp}
Let $n\in\omega$. A space $X$ is a subspace of an Eberlein compact (respectively, CC; UEC; $n$-UEC)  iff it is an E-space (respectively, C-space; UE-space; $n$-UE-space).
 \end{cor}

 \doc It follows from Theorem \ref{MainE} and the fact that the closed subspaces of Eberlein compacts (respectively, CC; UEC; $n$-UEC) are Eberlein compacts (respectively, CC; UEC; $n$-UEC).
 \sqs

\begin{cor}\label{metp}{\rm (\cite{P})}
 Every metrizable space has a compactification which is a 0-uniform Eberlein compact.
 \end{cor}

 \doc It follows from Theorem \ref{MainE},  the Bing Metrization Theorem and the fact that in metrizable spaces all open subsets are cozero-sets.
 \sqs

\begin{rem}\label{baleir}
\rm
 In  \cite{D1}, I mentioned that Theorem \ref{MainE} implies the Arhangel'ski\v{\i} Theorem \cite[Theorem 14]{A} but I missed to note that Theorem \ref{MainE} implies even a stronger result, namely that every metrizable space has a compactification which is a uniform Eberlein compact (the proof of both these results coincide with the proof of Corollary \ref{metp}).   I discovered this corollary a little bit  later  (\cite{D0}) and I intended to include it in the full version of the paper \cite{D1}. The writing of the full version was always postponed and in the meantime this assertion was noted by T. Banakh and A. Leiderman \cite{BL} (and supplied with a short analytic proof). Finally, B. A. Pasynkov \cite{P}  obtained the result stated in  Corollary \ref{metp}, which is stronger than both previous assertions.
In fact,  B. A. Pasynkov  proved in \cite{P} a theorem which is even  stronger than Corollary \ref{metp} (see \cite[Theorem 4]{P}).
 \end{rem}

 \begin{cor}\label{be}{\rm (\cite{D1})}
Every Baire E-space contains a dense $G_\d$ metrizable subspace.
 \end{cor}

 \doc It follows from theorems \ref{MainE} and \ref{Baire}.
 \sqs

 \begin{fact}\label{unifthSEC}
A compact Hausdorff space $X$ is a strong Eberlein compact iff it has a strongly point-finite family  $\a$ consisting of clopen sets such that $\a\cup\{X\stm V\st V\in\a\}$ is a subbase of $X$.
\end{fact}

 \doc  ($\Rightarrow$) Let $X$ be a SEC. Then $X$ has a $T_0$-separating, point-finite family $\a$ of clopen sets. Since $\a$ consists of closed subsets of $X$ and $X$ is compact, we get that $\a$ is  a strongly point-finite family. Using again the compactness of $X$, we obtain that the family $\a\cup\{X\stm V\st V\in\a\}$ is a subbase of $X$.

 \medskip

 \noindent($\Leftarrow$) Let $\a$ be a strongly point-finite family of clopen subsets of $X$ such that $\a\cup\{X\stm V\st V\in\a\}$ is a subbase of $X$.
Then, clearly, $\a$ is a $T_0$-separating family. Therefore, $X$ is a SEC.
\sqs

\begin{theorem}\label{MainSE}{\rm (\cite{D1})}
 A space X has a compactification $cX$ which is a strong
Eberlein compact iff it is an SE-space.
\end{theorem}

\doc ($\Rightarrow$)  Let $X$ has a compactification $cX$ which is a strong Eberlein compact. We can think that $X$ is a subspace of $cX$.
By Fact \ref{unifthSEC}, $cX$ has a a strongly point-finite family  $\a\ap$ consisting of clopen subsets of $cX$ such that $\a\ap\cup\{cX\stm V\st V\in\a\ap\}$ is a subbase of $cX$. Then it is easy to see that $\a=\a\ap\cap X$ is a strongly point-finite family consisting of clopen subsets of $X$ such that $\a\cup\{X\stm V\st V\in\a\}$ is a subbase of $X$.
Hence, $X$ is an SE-space.

 \medskip

 \noindent($\Leftarrow$) Let $X$ be an SE-space and $\a$ be a strongly point-finite family of clopen subsets of $X$ such that $\a\cup \{X\stm V\st V\in\a\}$ is a subbase of $X$.
Setting, for every $V\in\a$ and every $n\in\NNNN$,  $U_n(V)=V$,  we get, obviously,   a U-representation $U(\a)$  of $\a$ which certify that
 $\a$ is an almost subbase of $X$.
Let, for every $V\in\a$, $f_V:X\lra\DDDD$ be the characteristic function of $V$. Note that the function $f_V$ corresponds to the U-representation $U(V)$ of $V$ (see Lemma \ref{Uf}).
Set $Y=\DDDD^\a$ and let $f:X\lra Y$ be the diagonal of the family  $\{f_V\st V\in\a\}$.
Now, exactly as in the proof of Theorem \ref{MainE}, we show that  $f:X\lra Y$ is an embedding. Set
  $cX=\cl_Y(f(X))$.  We will show that $cX$ is a strong Eberlein compact. As in the proof of Theorem \ref{MainE}, we define the functions $f_V\ap:cX\lra\DDDD$ for every $V\in\a$, and then, using Lemma \ref{Ex},  we obtain that the family $\a\ap=\{coz(f_V\ap)\st V\in\a\}$ is strongly point-finite. Since, for every $V\in\a$, $coz(f_V\ap)=(f_V\ap)\inv(1)$, we get that $\a\ap$ consists of clopen subsets of $cX$. Also, it is clear that $\a\ap$ is a $T_0$-separating family. Therefore, $cX$
is a strong Eberlein compact.
 \sqs

\begin{rem}\label{reme1}
\rm
In connection with the next corollary, let us note that, obviously, the closed subspaces of $n$-UECs are $n$-UECs (for every $n\in\omega$); also, the countable product of $n$-UECs is an $n$-UEC (for every $n\in\omega$). The proof of the last assertion is similar to the proof of Proposition 2 of \cite{P}. Let us sketch it. Let $n\in\omega$ and $X=\prod\{X_i\st i\in\NNNN\}$, where, for every $i\in\NNNN$, $X_i$ is an $n$-UEC. Then, for every $i\in\NNNN$, there exists a $\s$-(n+1)-boundedly point-finite family $\g_i$ of cozero subsets of $X_i$ which $T_0$-separates $X_i$. Let, for every $i\in\NNNN$, $\pi_i:X\lra X_i$ be the projection. Set, for every $i\in\NNNN$, $\d_i=\{\pi_i\inv(U)\st U\in\g_i\}$ and $\d=\bigcup\{\d_i\st i\in\NNNN\}$. Then $\d$ is $\s$-(n+1)-boundedly point-finite family  of cozero subsets of $X$ which $T_0$-separates $X$. Thus, $X$ is an $n$-UEC.
\end{rem}

\begin{cor}\label{e1}{\rm (\cite{D1})} Let $n\in\omega$. Then:

\smallskip

\noindent(a) The property of being an E-space (respectively, C-space, UE-space, $n$-UE-space, SE-space) is
hereditary and additive;

\smallskip

\noindent(b) The property of being E-space (respectively, C-space, UE-space, $n$-UE-space) is
$\aleph_0$-multiplicative;

\smallskip

\noindent(c) If X is a C-space (respectively, E-space, UE-space, $n$-UE-space, SE-space), then $X$ is a Fr\'{e}chet-Urysohn space;

\smallskip

\noindent(d) If X is an E-space, then  $c(X)=w(X)=d(X)$, $c(X) \le | X | \le (c(X))^{\aleph_0}$,
$\aleph_0\le \chi(X)\le c(X)$ and the bounds given in the last two inequalities are the best possible;

\smallskip

\noindent(e) The following conditions are equivalent for a space $X$: (i) $X$ is separable metrizable; (ii) $X$ is a  0-UE-space with $c(X)=\aleph_0$; (iii) $X$ is a UE-space with $c(X)=\aleph_0$;  (iv) $X$
is an E-space with $c(X)=\aleph_0$.
 \end{cor}

\doc (a) The fact that the corresponding property is  hereditary is obvious. Let $J$ be a set, $\{X_j\st j\in J\}$ be a disjoint family of E-spaces and $X=\bigoplus\{X_j\st j\in J\}$. Let, for every $j\in J$, $\a_j$ be a $\s$-strongly point-finite almost subbase of $X_j$. Then $\a=\{X_j\st j\in J\}\cup\bigcup\{\a_j\st j\in J\}$ is a $\s$-strongly point-finite almost subbase of $X$. Hence, $X$ is an E-space. The proof for the other three cases is analogous.

\smallskip

\noindent(b) It follows from Theorem \ref{MainE}  and the fact that the class of ECs (respectively, CCs, UECs, $n$-UECs) is closed under countable products (see \cite[Proposition 3.3]{L}, \cite[Theorem 3.6]{W} and Remark \ref{reme1}).

\smallskip

\noindent(c) It follows from Theorem \ref{MainE}, the fact that Corson compacts are Fr\'{e}chet-Urysohn spaces (see, e.g., \cite[3.10.D]{E}) and \cite[2.1.H(b)]{E}.

\smallskip

\noindent(d) By Theorem \ref{MainE}, $X$ has a compactification $Y$ which is an EC. Then $c(Y)=w(Y)$ (see \cite[Theorem 3.1]{W}). Using \cite[2.7.9(d)]{E}, we get that $w(X)\le w(Y)=c(Y)=c(X)\le w(X)$. Thus $c(X)=w(X)$. Since $c(X)\le d(X)\le w(X)$, we get that $d(X)=w(X)=c(X)$. Further, by \cite[Theorem 3.1]{W}, $c(Y)\le |Y|\le (c(Y))^{\aleph_0}$. Then $c(X)\le |X|\le |Y|\le(c(Y))^{\aleph_0}=(c(X))^{\aleph_0}$. Finally, we have that $\aleph_0\le\chi(X)\le w(X)=c(X)$. The rest follows from
 the remark after \cite[Theorem 3.1]{W}.
\smallskip

\noindent(e) It follows from (d) and Corollary \ref{metp}.
\sqs

\begin{cor}\label{e2}
 Let $n\in\omega$ and $X$ be a locally compact Hausdorff space. Then:

 \smallskip

\noindent(a){\rm (\cite{D1})} the Alexandroff
compactilication $\a X$ of $X$ is an EC (respectively, CC, UEC, SEC) iff $X$ is an E-space (respectively, C-space, UE-space, SE-space);

\smallskip

\noindent(b) if, in addition, $X$ is a paracompact space   then $\a X$ is an $n$-UEC iff $X$ is an $n$-UE-space; thus, if $X$ is metrizable, then $\a X$ is a 0-UEC;

\smallskip

\noindent(c){\rm (\cite{D1})} if $Y$ is a perfect image of $X$ and $X$ is an E-space
(respectively, C-space, UE-space, SE-space) then $Y$ is an E-space (respectively, C-space, UE-space, SE-space).
\end{cor}

\doc (a) The necessity follows from theorems \ref{MainE} and \ref{MainSE}.
For the sufficiency, let $X$ be an E-space. Then, by Theorem \ref{MainE}, $X$ has a compactification $cX$ which is an EC. Since $\a X$ is a continuous image of $cX$,  Theorem \ref{contim} implies that $\a X$ is an EC. The proof for the other three cases is analogous.

\smallskip

\noindent(b) The necessity is clear. Let us prove the sufficiency. Let $X$ be a paracompact locally compact Hausdorff $n$-UE-space. By the Morita Theorem (see \cite[Theorem 5.1.27]{E}), we have that $X=\bigoplus\{X_j\st j\in J\}$, where $J$ is some set and, for every $j\in J$, $X_j$ is a Lindel\"{o}f space. Then, for every $j\in J$, every open $F_\s$-subset of $X_j$ is a union of countably many compact subsets of $X_j$, and thus it is a cozero-subset of $\a X$. Therefore, if, for every $j\in J$, $\g_j$ is a $\s$-(n+1)-boundedly point-finite family of cozero-subsets of $X_j$ that $T_0$-separates $X_j$, then $\g=\{X_j\st j\in J\}\cup\bigcup\{\g_j\st j\in J\}$ is a $\s$-(n+1)-boundedly point-finite family of cozero-subsets of $\a X$ that $T_0$-separates $\a X$. Hence, $\a X$ is an $n$-UEC.

Since every metrizable space $X$ is a 0-UE-space (see Corollary \ref{metp}),  the above result implies that if, in addition, $X$ is locally compact then $\a X$ is a 0-UEC.

\smallskip

\noindent(c) Let $f:X\lra Y$ be a perfect surjection. Then $f$ can be extended to a continuous map $f^\a:\a X\lra\a Y$ and all follows from (a) and Theorem \ref{contim}.
\sqs

\begin{rem}\label{reme2}
\rm
(a) Note that in the proof of Corollary \ref{e2}(a) we used Theorem \ref{contim} which is a very deep result. If we suppose in Corollary \ref{e2}(a) that $X$ is, in addition, a paracompact space, then the obtained weaker assertion for E-spaces, C-spaces and UE-spaces can be proved as Corollary \ref{e2}(b).

\smallskip

\noindent(b)
If the continuous image of an $n$-UEC, where $n\in\omega$, were an $n$-UEC, then we could remove the requirement of paracompactness in Corollary \ref{e2}(b) arguing as in the proof of Corollary \ref{e2}(a). But, {\em in general, the continuous image of an $n$-UEC, where $n\in\omega$, is not an $n$-UEC.} Indeed, let, for every cardinal number $\tau\ge\aleph_0$, $D_\tau$ denote the discrete space of cardinality $\tau$. Then $\a D_\tau$ is a 0-UEC. Hence, by Remark \ref{reme1}, $X_\tau=(\a D_\tau)^{\aleph_0}$ is a 0-UEC, as well as every closed subset of $X_\tau$. By \cite[Lemma 1.2]{BRW}, for every UEC $Y$ there exists a $\tau\ge\aleph_0$ and a closed subset $F$ of $X_\tau$ such that $Y$ is a continuous image of $F$. Thus, {\em every UEC is a continuous image of a 0-UEC}. Since there exist UECs which are not 0-UECs (see \cite{P}), we obtain that the continuous image of a 0-UEC is not, in general, a 0-UEC.
\end{rem}

\begin{cor}\label{e3}{\rm (\cite{L})}
If $X$ is locally compact metrizable or if $X$ is a disjoint sum of ECs
then the Alexandroff compactilication $\a X$ of $X$ is an EC.
\end{cor}

\doc It follows from corollaries
\ref{metp},
\ref{e1}(a) and Remark \ref{reme2}(a).
\sqs

\begin{cor}\label{e30}
Let $n\in\omega$. If $X$ is  a disjoint sum of $n$-UECs
then the Alexandroff compactilication $\a X$ of $X$ is an $n$-UEC.
\end{cor}

\doc It follows from  \cite[Theorem 5.1.30]{E} and corollaries
\ref{e1}(a), \ref{e2}(b).
\sqs

\begin{defi}\label{meta}{\rm (\cite{D1})}
\rm
 A space $X$ is called {\em $\s$-strongly metacompact}\/ (respectively, {\em strongly
metalindel\"{o}f})\/ if every open cover of $X$ has a $\s$-strongly point-finite (respectively, strongly
point-countable) open refinement.
\end{defi}

In what follows, if $\AA$ is a family of subsets of a space $X$ then we shall denote by $\overline{\AA}$ the family $\{\cl_X(A)\st A\in\AA\}$.

\begin{lm}\label{lmlc}{\rm (\cite{D0})}
 Let $X$ be a locally compact Hausdorff space. Then:

 \smallskip

\noindent(a) $X$  is
(hereditarily) strongly metalindel\"{o}f iff it is (hereditarily)  metalindel\"{o}f;

 \smallskip

\noindent(b) $X$ is (hereditarily) $\s$-strongly metacompact iff it is  (hereditarily)  $\s$-metacompact.
\end{lm}

\doc (a) The necessity is obvious. For the sufficiency, let us first regard the case when $X$ is a  metalindel\"{o}f locally compact Hausdorff space. Let $\UU$ be an open cover of $X$. Then $\UU$ has an open refinement $\UU\ap$ whose elements have compact closures. Let $\VV$ be an open point-countable refinement of $\UU\ap$. Then, by \cite[Theorem 1.1 and Footnote 1]{GM}, there exists an open cover $\WW=\{W_V\st V\in\VV\}$ of $X$ such that $\cl_X(W_V)\sbe V$ for every $V\in\VV$. Let $\mu\sbe\WW$ and $|\mu|=\aleph_1$. Then $\bigcap\overline{\mu}=\ems$ because $\overline{\WW}$ is a
shrinking of $\VV$ and $\VV$ is point-countable. Since the elements of $\overline{\mu}$ are compact sets, there exists a finite subfamily $\mu_0$ of $\mu$ such that $\bigcap\overline{\mu_0}=\ems$. Then $\bigcap\mu_0=\ems$. Hence, the cover $\WW$ is strongly point-countable. This implies that $X$ is strongly metalindel\"{o}f.

Let now $X$ be  a hereditarily metalindel\"{o}f locally compact Hausdorff space. Then, clearly, by the above paragraph, every open subspace of $X$ is strongly
metalindel\"{o}f. This implies easily that $X$ is hereditarily strongly metalindel\"{o}f.

 \smallskip

\noindent(b)  The proof is similar to that one of (a).
\sqs

\begin{rem}\label{remp}
\rm
A space $X$ is called {\em strongly metacompact}\/ (\cite{Pe1}) (note that  the name $``$b-paracompact space" is used in \cite{Pe1}) if
every open cover of $X$ has a strongly point-finite  open refinement. In \cite{Pe1}, it is proved that if $X$ is a locally compact Hausdorff space, then
$X$ is  strongly metacompact iff it is    metacompact (as it is written there, this is a result of M. Patashnik). Obviously, this implies that if $X$ is a locally compact Hausdorff space, then
$X$ is hereditarily strongly metacompact iff it is  hereditarily  metacompact.

Let us also note that in \cite{Po} the results stated in  Lemma \ref{lmlc} were mentioned using \cite{D0} as a prime source.
\end{rem}

Using Lemma \ref{lmlc} and Yakovlev's Theorem (\cite[Corollary 2]{Y}), we obtain:

\begin{cor}\label{e4}{\rm (\cite{D1})}
 Every Eberlein space (respectively, C-space) is hereditarily $\s$-strong\-ly me\-tacompact
(respectively, hereditarily strong\-ly metalindel\"{o}f).
\end{cor}

\begin{rem}\label{remstr}
\rm
Recall that, by the Gruenhage Theorem \cite[Theorem 2.2]{Gr}, a compact Hausdorff space is an Eberlein compact iff $X^2$ is hereditarily $\s$-metacompact. In 1991, I asked  S. Popvassilev (who  was my student at that time) {\em whether a Tychonoff space $X$ is an E-space iff $X^2$ is hereditarily $\s$-strongly metacompact}; in his Master Thesis (University of Sofia, 1992) he observed that {\em the space $X=C_p(\IIII)$ is  a non-E-space}\/ (because $c(X)=\aleph_0<2^{\aleph_0}=|\IIII|=w(X)$ by \cite[Theorems 0.3.7 and I.1.1]{A2}, and this contradicts  Corollary \ref{e1}(d)) {\em but $X^2$ is hereditarily $\s$-strongly metacompact}\/ (indeed, $X^2$ is  hereditarily Lindel\"{o}f by the Zenor-Velichko Theorem (see \cite[Theorem II.5.10]{A2}); thus, by \cite[Corollary 5.3.11]{E}, $X^2$ is hereditarily strongly paracompact and, hence, $X^2$ is hereditarily ($\s$-)strongly metacompact).
\end{rem}

S. A. Peregudov  \cite[Lemma 3]{Pe1}  proved that in every space $X$ the cardinality of every strongly point-finite family of open sets is less or equal to $c(X)$.
 This result implies that if $X$ is $\s$-strongly metacompact then $l(X)\le c(X)$  (see \cite{Po}) (here, as usual, $l(X)$ is the Linel\"{o}f number of $X$,
i.e., the smallest cardinal number $\tau$ such that every open cover of $X$ has an open refinement of cardinality $\le\tau$). Thus, using Corollary \ref{e4}, we obtain the following fact:

\begin{cor}\label{e5}
 If $X$ is an  E-space, then $l(X)\le c(X)$.
\end{cor}

Let us note that S. A. Peregudov \cite[Theorem 2]{Pe1} proved that if $X$ is a strongly metacompact space then $l(X)\le c(X)$.

\begin{cor}\label{e6}
 If $X$ is an  E-space, then $hl(X)= c(X)=hc(X)$.
\end{cor}

\doc By Corollary \ref{e1}(d), we have that $c(X)=w(X)$. Since, obviously, $c(X)\le hc(X)\le w(X)$, we get that $c(X)=hc(X)$. Let $A\sbe X$. Then, by Corollaries \ref{e1}(a) and \ref{e5}, $l(A)\le c(A)\le hc(X)=c(X)$. Therefore, $hl(X)\le c(X)$. Since $hl(X)\ge c(X)$ (see, e.g., \cite[3.12.7(e)]{E}), we get that
 $hl(X)= c(X)$.
 \sqs

\begin{rem}\label{remunif}
\rm
By the Arhangel'ski\v{\i} Theorem \cite[Theorem 14]{A}, every metrizable space is an E-space.
On the other hand, as it is noted in \cite{A4}, $``$not
every completely regular Moore space has an Eberlein compactification, since
there are separable non-metrizable Moore spaces". Such a space is, for example, the Niemytzki plane (see, e.g., \cite{Bu}). By the Yakovlev Theorem (\cite[Corollary 2]{Y}), every Eberlein compact is hereditarily $\s$-metacompact.  Hence, it is natural to ask whether a Tychonoff (or normal)  metacompact Moore space
is an E-space. Since a regular space has a uniform base (in the sense of P. S. Alexandrov \cite{Alex}) iff it is a metacompact Moore space (\cite[Corollary and Theorem III]{Alex}, \cite[Theorem 4]{He}), this question is equivalent to the following one: {\em is every Tychonoff (or normal  $T_1$) space with a uniform base
 an E-space?}
 (Let us recall that, by \cite{Alex},  a space is metrizable iff it is a collectionwise normal $T_1$-space with a uniform base.)

The Niemytzki plane is not a space with a uniform base  because it is not a metacompact space (see \cite[Exercise 5.3.B(a)]{E}).
Let us show that
   the Pixley-Roy example $X$  of a Tychonoff non-separable Moore metacompact space with $c(X)=\aleph_0$ (see \cite{PR}) {\em is a Tychonoff non-E-space with a uniform base.}  Indeed, we have that: (a) $X$ is a space with a uniform base,\\
    (b) $X$ is a p-space (in the sense of \cite{A3}) (because each completely regular Moore space is a p-space (\cite{A1} or \cite[V.226]{AP})),\\
     (c) as a space with a uniform base, $X$ has a point-countable base (see  \cite[Prop. I]{Alex}).\\
    Supposing  that $X$ is a Lindel\"{o}f space, we obtain that $X$ is metrizable (because each paracompact p-space with a point-countable base is metrizable (see \cite{Pon} or \cite[V.229]{AP})), and, thus,
 $c(X)=d(X)$, a contradiction. Hence, $X$ is a non-Lindel\"{o}f space, i.e., $l(X)>\aleph_0$. Since $c(X)=\aleph_0$,  Corollary \ref{e5} implies that $X$ is a non-E-space.

 Starting with the Przymusi\'{n}ski-Tall example of a normal  non-separable Moore metacompact $T_1$-space $X$ with $c(X)=\aleph_0$, constructed under  MA+$\neg$CH in \cite{PT}, we obtain, as above, that $X$ is a non-E-space. Therefore,
 {\em under MA+$\neg$CH,  there exists a normal Hausdorff non-E-space with a uniform base.}

P. J. Nyikos \cite{Ny} proved that the Product Measure Extension
Axiom (PMEA, for short) implies that  normal Moore $T_1$-spaces are metrizable. Thus, {\em under PMEA,
every normal $T_1$-space with a uniform base}\/ is metrizable and, hence, it {\em is an E-space.}
\end{rem}

\section{On spaces co-absolute with Eberlein compacts and Ponomarev's solution of Birkhoff's Problem 72}

The fact that two Boolean algebras $A$ and $B$ are isomorphic  will be expressed by $``A\cong B$".

\begin{nota}\label{bir}
\rm We shall denote by:
\begin{itemize}
\item  $\MM$
the class of all metrizable
spaces,
\item  $\EE$
the class of all Eberlein compacts.
\end{itemize}

If $\KK$ is a class of topological spaces, we will set
$$\BB\KK=\{A\in |\Bool|\st (\ex X\in\KK)(A\cong RO(X)) \}.$$
\end{nota}

Recall that  a subset $B$ of a Boolean algebra $A$ is said to be {\em $\s$-disjointed}\/ (respectively, {\em dense}\/) if   $B=\bigcup\{B_n\st
n\in\mathbb{N}\}$, where for every $n\in\mathbb{N}$ and for
every two different elements $a,b$ of $B_n$ we have $a\we b=0$ (respectively, if for any $a\in A\stm\{0\}$
there exists $b\in B\stm\{0\}$ such that $b\le a$).

The Problem 72 of G. Birkhoff \cite{Bi} is the following:
characterize internally the elements of the class $\BB\MM$. It was
solved by V. I. Ponomarev \cite{P1}. He proved the following
beautiful theorem: if $A$ is a complete Boolean algebra then
$A\in\BB\MM$ iff it has a $\s$-disjointed dense subset $B$.
The proof of this theorem is difficult. We will obtain a  shorter
 proof of it on the base of some results from \cite{A,Nam,D2} which appeared many years after the publication of
Ponomarev's paper \cite{P1}.

 We will need a lemma from \cite{CNG}:

\begin{lm}\label{isombool}
Let $X$ be a dense subspace of a topological space $Y$. Then the
 Boolean
algebras $RC(X)$ and $RC(Y)$ are isomorphic.
\end{lm}

If $X$ is a set and $\g=\bigcup\{\g_n\st n\in\NNNN\}$, where, for every $n\in\NNNN$, $\g_n$ is a disjoint family of subsets of $X$, then we will say that $\g$ is a
{\em $\s$-disjoint family}. (Obviously, $``\s$-disjoint family" = $``\s$-1-boundedly point-finite family".)

The new proof of the Ponomarev Theorem as well as the  characterization of the
class of spaces which are co-absolute with (zero-dimensional)
Eberlein compacts will follow from the next theorem.

\begin{theorem}\label{birkhoffpon}
 A complete Boolean algebra $A$ is isomorphic to an algebra
of the form $RC(X)$, where $X$ is a (zero-dimensional) Eberlein
compact, iff $A$ has a $\s$-disjointed dense subset.
\end{theorem}

\doc ($\Rightarrow$) Let $A$ be a Boolean algebra which is
isomorphic to RC(X), where $X$ is an  Eberlein compact. By Corollary \ref{gdelta}, there exists a metrizable dense subset $Y$ of
$X$. Hence, by Lemma \ref{isombool}, $A$ is isomorphic to $RC(Y)$. The space $Y$ has a
$\sigma$-discrete base $\BB=\bigcup\{\BB_i\st i\in\mathbb{N}^+\}$,
where $\BB_i$ is a discrete family for every $i\in\mathbb{N}^+$.
Set, for every $i\in\mathbb{N}^+$, $\BB_i\ap=\{\cl(U)\st
U\in\BB_i\}$, and let $\BB\ap=\bigcup\{\BB_i\ap\st
i\in\mathbb{N}^+\}$. Then, obviously, $\BB\ap$ is a
$\sigma$-disjointed dense subset of $RC(Y)$. Hence, $A$  has a
$\s$-disjointed dense subset.

\smallskip

\noindent($\Leftarrow$) Let $A$ be a complete Boolean algebra
having a $\s$-disjointed dense subset $B_0$. Let $B$ be the
Boolean subalgebra of $A$ generated by $B_0$. Then $A$ is a
minimal completion of $B$. Set $X=S(B)$.
 Then $X$ is a
zero-dimensional compact Hausdorff space and there exists an
isomorphism $\p:B\lra CO(X)$. We will show that $\BB=\p(B_0)$ is a
$\s$-disjoint almost subbase of $X$.  For every $V\in\BB$ and
every $n\in\mathbb{N}^+$, set $U_n(V)=V$. Then $\{U_n(V)\st
i\in\mathbb{N}^+\}$ is a U-representation of $V$. Hence, it is enough
 to show that the family $\BB\ap=\BB\cup\{X\stm V\st
V\in\BB\}$ is a subbase of $X$. Obviously, $\BB\ap=\p(B_0\cup
B_0^*)$, where $B_0^*=\{b^*\st b\in B_0\}$. Since the
set of all finite joins of all finite meets of the elements of the
subset $B_0\cup B_0^*$ of $A$ coincides with $B$ (see \cite[Proposition 4.4]{kop89}),
we get that the
family of all finite unions of the finite intersections of the
elements of the family $\BB\ap$ coincides with $CO(X)$. The family $CO(X)$ is a
base of $X$; hence, the family of all finite intersections of the
elements of $\BB\ap$ is a base of $X$, i.e., $\BB\ap$ is a subbase
of $X$. Therefore, $\BB$ is an almost subbase of $X$. Since $\BB$
is, obviously, a $\s$-disjoint family, we get, by Theorem
\ref{unifth} (or by Theorem
\ref{unifthE}), that $X$ is an Eberlein compact. Now, $RC(X)$ is a
minimal completion of $CO(X)$; thus  $RC(X)$ and $A$ are
isomorphic Boolean algebras. \sqs

\begin{rem}\label{rempp}
\rm
Note that the zero-dimensional Eberlein compact $X$ constructed in the proof of the sufficiency of the preceding theorem is even a 0-UEC (by Theorem \ref{unifthE}).
\end{rem}

\begin{pro}\label{birk0}
$\BB\MM=\BB\EE$.
\end{pro}

\doc By the  Arhangel'ski\v{i} Theorem \cite[Theorem 14]{A},
 every metric space can be
densely embedded in an Eberlein compact. Conversely,
by Corollary \ref{gdelta},
every Eberlein compact contains a dense  metrizable subspace.
Applying  Lemma \ref{isombool}, we conclude that $\BB\MM=\BB\EE$.
\sqs

Combining  Theorem \ref{birkhoffpon} with Proposition \ref{birk0}, we obtain
the Ponomarev Theorem \cite{P1} giving a solution of Birkhoff's
Problem 72 \cite{Bi}.

\begin{cor}\label{birkhoffponom}{\rm (V. I. Ponomarev \cite{P1})}
 A complete Boolean algebra $A$ is isomorphic to an algebra
of the form $RC(X)$, where $X$ is a metrizable space, iff $A$ has
a $\s$-disjointed dense subset.
\end{cor}

Recall  that if $X$ is a regular space then a space $EX$ is
called an {\em absolute} of $X$ iff there exists  a perfect
irreducible map $\pi_X:EX\lra X$ and every perfect  irreducible
 preimage of $EX$ is homeomorphic to $EX$ (see, e.g., \cite{PS}).
 Two regular spaces are said to be {\em co-absolute}\/ if their
 absolutes are homeomorphic.
 It
is well-known that:  a) the absolute is unique up to
homeomorphism; b) a space $Y$ is an absolute  of a regular space
$X$ iff $Y$ is an extremally disconnected Tychonoff space
 for which there exists a perfect irreducible map
$\pi_X:Y\lra X$;  c) if $X$ is a compact Hausdorff space then
$EX=S(RC(X))$ (see, e.g., \cite{Walker}).

\begin{theorem}\label{abseb}
Let $X$ be a compact Hausdorff space. Then the following
conditions are equivalent:

\smallskip

\noindent(a) $X$ is co-absolute with an Eberlein compact;

\smallskip

\noindent(b) $X$ has a $\s$-disjoint $\pi$-base;

\smallskip

\noindent(c)  $X$ is co-absolute with a zero-dimensional 0-UEC;

\smallskip

\noindent(d)  $X$ is co-absolute with a zero-dimensional Eberlein
compact.
\end{theorem}

\doc (a)$\Rightarrow$(b) Let $Y$ be an Eberlein compact which is co-absolute
with $X$. Then $RC(Y)\cong RC(X)$. By Theorem \ref{birkhoffpon}, we get that
the Boolean algebra $RO(X)$ (which is isomorphic to the Boolean
algebra $RC(X)$) has a $\s$-disjointed dense subset $\AA$. Then,
obviously, $\AA$ is a $\s$-disjoint $\pi$-base of $X$.

\smallskip

\noindent(b)$\Rightarrow$(c) Let $\AA$ be a  $\s$-disjoint $\pi$-base of
$X$. Set $\AA\ap=\{\int(\cl(U))\st U\in\AA\}$. Then, obviously,
$\AA\ap$ is a $\s$-disjointed dense subset of the Boolean algebra
$RO(X)$. Since $RO(X)\cong RC(X)$, Remark \ref{rempp}
implies that there exists a zero-dimensional 0-UEC $Y$
with $RC(Y)\cong RC(X)$. Thus $X$ and $Y$ are co-absolute spaces.

\smallskip

\noindent(c)$\Rightarrow$(d) and (d)$\Rightarrow$(a) Both implications are obvious. \sqs


\baselineskip = 0.75\normalbaselineskip

\end{document}